\documentclass[12pt,a4paper]{article}
\usepackage{amsmath}
\usepackage{amssymb}
\usepackage{t1enc}
\usepackage[latin1]{inputenc}
\usepackage[english]{babel}
\usepackage{amsfonts}
\usepackage{latexsym}
\usepackage{bm}
\usepackage{setspace}
\usepackage[toc,page]{appendix}
\usepackage{graphicx}
\usepackage{enumerate}
\usepackage[numbers,sort&compress]{natbib}
\usepackage{tikz}

\usepackage{theoremref}
\usepackage{lineno}

\usepackage{mathtools}

\newtheorem{theorem}{Theorem}
\newtheorem{prop}{Proposition}
\newtheorem{lemma}{Lemma}

\newtheorem{conj}{Conjecture}

\newtheorem{construction}{Construction}
\newtheorem{problem}{Problem}

\allowdisplaybreaks

\begin{document}
\title{Isolation of regular graphs, stars and $k$-chromatic graphs \medskip\medskip
}

\author{Peter Borg \\[2mm]
\normalsize Department of Mathematics \\
\normalsize Faculty of Science \\
\normalsize University of Malta\\
\normalsize Malta\\
\normalsize \texttt{peter.borg@um.edu.mt}
}

\date{}
\maketitle

\begin{abstract}
Given a set $\mathcal{F}$ of graphs, we call a copy of a graph in $\mathcal{F}$ an $\mathcal{F}$-graph. The $\mathcal{F}$-isolation number of a graph $G$, denoted by $\iota(G,\mathcal{F})$, is the size of a smallest set $D$ of vertices of $G$ such that the closed neighbourhood of $D$ intersects the vertex sets of the $\mathcal{F}$-graphs contained by $G$ (equivalently, $G - N[D]$ contains no $\mathcal{F}$-graph). Thus, $\iota(G,\{K_1\})$ is the domination number of $G$, and $\iota(G,\{K_2\})$ is the vertex-edge domination number of $G$. Clearly, $\iota(G, \mathcal{F}) \leq \iota(G, \mathcal{F} \cup \mathcal{H})$. For any integer $k \geq 1$, let $\mathcal{F}_{0,k}$ be the set consisting of the $k$-star $K_{1,k}$, let $\mathcal{F}_{1,k}$ be the set of regular graphs whose degree is at least $k-1$, let $\mathcal{F}_{2,k}$ be the set of graphs whose chromatic number is at least $k$, and let $\mathcal{F}_{3,k}$ be the union $\mathcal{F}_{0,k} \cup \mathcal{F}_{1,k} \cup \mathcal{F}_{2,k}$. We prove that if $G$ is a connected $n$-vertex graph, then $\iota(G, \mathcal{F}_{3,k}) \leq \frac{n}{k+1}$ unless $G$ is a $k$-clique or $k = 2$ and $G$ is a $5$-cycle. This generalizes a classical bound of Ore on the domination number, a bound of Caro and Hansberg and of \.{Z}yli\'{n}ski on the vertex-edge domination number, a bound of Fenech, Kaemawichanurat and the author on the $k$-clique isolation number, a bound of the author on the cycle isolation number, and a bound of Caro and Hansberg on the $\mathcal{F}_{0,k}$-isolation number. The proof features a new strategy. For $i = 1, 2, 3$, the bound $\frac{n}{k+1}$ on $\iota(G, \mathcal{F}_{i,k})$ is attainable if $k+1$ divides $n$. Our second main result is that the bound $\frac{n}{k+1}$ on $\iota(G, \mathcal{F}_{0,k})$ is attainable if and only if $n$ is $0$ or $k+1$ or $2(k+1)$. We pose some problems and conjectures, and establish additional intriguing phenomena concerning $k$-star isolation and $k$-cycle isolation.
\end{abstract}

\section{Introduction} \label{Introsection}
For standard terminology in graph theory, we refer the reader to \cite{West}. Most of the notation and terminology used here is defined in \cite{Borg}. The set of positive integers is denoted by $\mathbb{N}$. For $n \in \{0\} \cup \mathbb{N}$, $[n]$ denotes the set $\{i \in \mathbb{N} \colon i \leq n\}$. Note that $[0]$ is the empty set $\emptyset$. Arbitrary sets and graphs are taken to be finite. For a set $X$, ${X \choose 2}$ denotes the set of $2$-element subsets of $X$. Every graph $G$ is taken to be \emph{simple}, that is, its vertex set $V(G)$ and edge set $E(G)$ satisfy $E(G) \subseteq {V(G) \choose 2}$. We may represent an edge $\{v,w\}$ by $vw$. We call $G$ an \emph{$n$-vertex graph} if $|V(G)| = n$. For a vertex $v$ of $G$, $N_{G}(v)$ denotes the set of neighbours of $v$ in $G$, $N_{G}[v]$ denotes the closed neighbourhood $N_{G}(v) \cup \{ v \}$ of $v$, and $d_{G}(v)$ denotes the degree $|N_{G} (v)|$ of $v$. For a subset $X$ of $V(G)$, $N_G[X]$ denotes the closed neighbourhood $\bigcup_{v \in X} N_G[v]$ of $X$, $G[X]$ denotes the subgraph of $G$ induced by $X$ (that is, $G[X] = (X,E(G) \cap {X \choose 2})$), and $G - X$ denotes the subgraph of $G$ obtained by deleting the vertices in $X$ from $G$ (that is, $G - X = G[V(G) \backslash X]$). We may abbreviate $G - \{x\}$ to $G - x$. Where no confusion arises, the subscript $G$ may be omitted from notation that uses it. The complete graph $([n], {[n] \choose 2})$ is denoted by $K_n$. 

If $G$ and $H$ are graphs, $f : V(H) \rightarrow V(G)$ is a bijection, and $E(G) = \{f(v)f(w) \colon vw$ $\in E(H)\}$, then we say that $G$ is a \emph{copy of $H$} or that $G$ is \emph{isomorphic to $H$}, and we write $G \simeq H$. Thus, a copy of $H$ is a graph obtained by relabelling the vertices of $H$. If $G$ and $H$ are graphs such that $V(H) \subseteq V(G)$ and $E(H) \subseteq E(G)$, then $H$ is called a \emph{subgraph of $G$}, and we say that \emph{$G$ contains $H$}.

If $D \subseteq V(G) = N[D]$, then $D$ is called a \emph{dominating set of $G$}. The size of a smallest dominating set of $G$ is called the \emph{domination number of $G$} and is denoted by $\gamma(G)$. If $\mathcal{F}$ is a set of graphs and $F$ is a copy of a graph in $\mathcal{F}$, then we call $F$ an \emph{$\mathcal{F}$-graph}. A subset $D$ of $V(G)$ is called an \emph{$\mathcal{F}$-isolating set of $G$} if $N[D]$ intersects the vertex sets of the $\mathcal{F}$-graphs contained by $G$. Thus, $D$ is an $\mathcal{F}$-isolating set of $G$ if and only if $G - N[D]$ contains no $\mathcal{F}$-graph. It is to be assumed that $(\emptyset, \emptyset) \notin \mathcal{F}$. The size of a smallest $\mathcal{F}$-isolating set of $G$ is called the \emph{$\mathcal{F}$-isolation number of $G$} and is denoted by $\iota(G, \mathcal{F})$. If $\mathcal{F} = \{F\}$, then we may replace $\mathcal{F}$ in these defined terms and notation by $F$. Clearly, $D$ is a dominating set of $G$ if and only if $D$ is a $K_1$-isolating set of $G$. Thus, $\gamma(G) = \iota(G, K_1)$.

The study of isolating sets was introduced by Caro and Hansberg~\cite{CaHa17}. It is an appealing and natural generalization of the classical domination problem \cite{C, CH, HHS, HHS2, HL, HL2}. One of the earliest results in this field is the upper bound $n/2$ of Ore \cite{Ore} on the domination number of any connected $n$-vertex graph $G$ with $n \geq 2$ (see \cite{HHS}). While deleting the closed neighbourhood of a dominating set yields the graph with no vertices, deleting the closed neighbourhood of a $K_2$-isolating set yields a graph with no edges. In the literature, a $K_2$-isolating set is also called a \emph{vertex-edge dominating set}. Caro and Hansberg~\cite{CaHa17} proved that if $G$ is a connected $n$-vertex graph with $n \geq 3$, then $\iota(G, K_2) \leq n/3$ unless $G$ is a $5$-cycle. This was independently proved by \.{Z}yli\'{n}ski \cite{Z} and solved a problem in \cite{BCHH}. The graphs attaining the bound have recently been partially determined by Lema\'{n}ska, Mora and Souto-Salorio \cite{LMS}, and subsequently fully determined by Boyer and Goddard \cite{BG}. Fenech, Kaemawichanurat and the present author~\cite{BFK} established the general sharp bound on $\iota(G, K_k)$ in Theorem~\ref{result1}, and this solved a problem of Caro and Hansberg~\cite{CaHa17}. In this paper, we show that sharp upper bounds on $\iota(G, \mathcal{F})$ that have been established for certain sets $\mathcal{F}$ hold for a significantly larger set of graphs. The bounds are sharp by Construction~\ref{const1}, given in \cite{BFK} for Theorem~\ref{result1}. 

For $n \geq 1$, the graphs $([n], \{\{1, i\} \colon i \in [n] \backslash \{1\}\})$ and $([n], \{\{i,i+1\} \colon i \in [n-1]\})$ are denoted by $K_{1,n-1}$ and $P_n$, respectively. For $n \geq 3$, $C_n$ denotes the graph $([n], \{\{1,2\}, \{2,3\}, \dots, \{n-1,n\}, \{n,1\}\})$. A copy of $K_n$ is called an \emph{$n$-clique} or simply a \emph{clique}. A copy of $K_{1,n}$ is called an \emph{$n$-star} or simply a \emph{star}. A copy of $P_n$ is called an \emph{$n$-path} or simply a \emph{path}. A copy of $C_n$ is called an \emph{$n$-cycle} or simply a \emph{cycle}. A $3$-cycle is a $3$-clique and is also called a \emph{triangle}. 

\begin{construction}[\cite{BFK}] \label{const1} \emph{Consider any $n, k \in \mathbb{N}$. Let $q = \left \lfloor \frac{n}{k+1} \right \rfloor$. Thus, $n = q(k+1) + r$, where $0 \leq r \leq k$. If $q = 0$, then let $B_{n,k} = K_r$. Suppose $q \geq 1$. Let $B_1, \dots, B_q$ be $(k+1)$-cliques, and let $B_{q+1}$ be an $r$-clique, such that $B_1, \dots, B_q, B_{q+1}$ are pairwise vertex-disjoint. Let $q' = q$ if $r = 0$, and let $q' = q+1$ if $r \geq 1$. For each $i \in [q']$, let $b_i \in V(B_i)$. Let $B_{n,k}$ be the connected $n$-vertex graph with $V(B_{n,k}) = \bigcup_{i=1}^{q'} V(B_i)$ and $E(B_{n,k}) = \{b_ib_{i+1} \colon i \in [q'-1]\} \cup \bigcup_{i=1}^{q'} E(B_i)$.}
\end{construction}

If $G$ is a graph, $k \geq 1$, and either $G$ is a $k$-clique or $k = 2$ and $G$ is a $5$-cycle, then we say that the pair $(G,k)$ is \emph{special}.

\begin{theorem}[\cite{BFK}] \label{result1}
If $G$ is a connected $n$-vertex graph and $(G,k)$ is not special, then
\[ \iota(G, K_k) \leq \left\lfloor \frac{n}{k+1} \right\rfloor. \]
Moreover, equality holds if $G = B_{n,k}$.
\end{theorem}
It is worth mentioning that the authors of \cite{BFK} also obtained a sharp upper bound in terms of the number of edges \cite{BFK2}. The above-mentioned result of Ore and result of Caro and Hansberg and of \.{Z}yli\'{n}ski are the cases $k = 1$ and $k = 2$ of Theorem~\ref{result1}, respectively.
 
Let $\mathcal{C}$ be the set of cycles. Solving another problem of Caro and Hansberg~\cite{CaHa17}, the present author~\cite{Borg} proved the following result.

\begin{theorem}[\cite{Borg}] \label{result2}
If $G$ is a connected $n$-vertex graph that is not a triangle, then
\[\iota(G,\mathcal{C}) \leq \left \lfloor \frac{n}{4} \right \rfloor.\] 
Moreover, equality holds if $G = B_{n,3}$.
\end{theorem}

The \emph{maximum degree of $G$}, denoted by $\Delta(G)$, is $\max\{d(v) \colon v \in V(G)\}$. For $k \geq 1$, a subset $D$ of $V(G)$ is a $K_{1,k}$-isolating set of $G$ if and only if $\Delta(G - N[D]) < k$. Caro and Hansberg \cite{CaHa17} proved the following result.

\begin{theorem}[\cite{CaHa17}] \label{result3} If $G$ is an $n$-vertex graph, then
\[\iota(G, K_{1,k}) \leq \left \lfloor \frac{n}{k+1} \right \rfloor.\]
Moreover, equality holds if $G$ is the union of pairwise vertex-disjoint copies of $K_{1,k}$.
\end{theorem}

Domination and isolation have been particularly investigated for maximal outerplanar graphs (mops) \cite{BK, BK2, CaWa13, CaHa17, Ch75, DoHaJo16, DoHaJo17, HeKa18, LeZuZy17, Li16, MaTa96, KaJi, To13}, mostly due to connections with Chv\'{a}tal's Art Gallery Theorem \cite{Ch75}. Kaemawichanurat and the present author~\cite{BK2} proved that $\iota(G, K_{1,k}) \leq \left \lfloor \frac{n}{k+3} \right \rfloor$ if $G$ is a mop.

If $V(G) \neq \emptyset$ and $d(v) = r$ for each $v \in V(G)$, then $G$ is said to be \emph{$r$-regular} or simply \emph{regular}, and $r$ is called the \emph{degree of $G$}. If there exists a function $f \colon V(G) \rightarrow [k]$ such that $f(v) \neq f(w)$ for every $v, w \in V(G)$ with $vw \in E(G)$, then $G$ is said to be \emph{$k$-colourable}. The smallest non-negative integer $k$ such that $G$ is $k$-colourable is called the \emph{chromatic number of $G$} and is denoted by $\chi(G)$. If $k = \chi(G)$, then $G$ is said to be \emph{$k$-chromatic}. Brooks' Theorem \cite{Brooks} tells us that for any connected $n$-vertex graph $G$, $\chi(G) \leq \Delta(G)$ unless $G$ is an $n$-clique or $n$ is odd and $G$ is a cycle.

For $k \geq 1$, let $\mathcal{F}_{0,k} = \{K_{1,k}\}$, let $\mathcal{F}_{1,k}$ be the set of regular graphs whose degree is at least $k-1$, let $\mathcal{F}_{2,k}$ be the set of graphs whose chromatic number is at least $k$, and let $\mathcal{F}_{3,k}$ be the union of $\mathcal{F}_{0,k}$, $\mathcal{F}_{1,k}$ and $\mathcal{F}_{2,k}$. In Section~\ref{Proofsection}, we prove the following result.

\begin{theorem} \label{mainresult1} If $G$ is a connected $n$-vertex graph and $(G,k)$ is not special, then
\[\iota(G,\mathcal{F}_{0,k} \cup \mathcal{F}_{1,k}) \leq \left \lfloor \frac{n}{k+1} \right \rfloor.\]
Moreover, equality holds if $G = B_{n,k}$.
\end{theorem}
Theorem~\ref{mainresult1} generalizes Theorems~\ref{result1}--\ref{result3} (as $k$-cliques are $(k-1)$-regular and cycles are $2$-regular), and the argument in its proof is more efficient than those in \cite{BFK, Borg} (for Theorems~\ref{result1} and \ref{result2}). In Section~\ref{Proofsection}, we also show that by Brooks' Theorem, Theorem~\ref{mainresult1} immediately yields the following generalization.

\begin{theorem} \label{mainresult1gen} If $i \in \{0, 1, 2, 3\}$, $G$ is a connected $n$-vertex graph and $(G,k)$ is not special, then 
\[\iota(G,\mathcal{F}_{i,k}) \leq \left \lfloor \frac{n}{k+1} \right \rfloor.\]
Moreover, equality holds if $i \in \{1, 2, 3\}$ and $G = B_{n,k}$.
\end{theorem}

In this line of study, a central aim is to determine if there exists a smallest constant $a(\mathcal{F})$ such that $a(\mathcal{F}) |V(G)|$ is an upper bound on $\iota(G,\mathcal{F})$ for every connected graph $G$ except for a finite number of non-isomorphic graphs and their copies. By Theorem~\ref{mainresult1gen}, $a(\mathcal{F}_{i,k}) = \frac{1}{k+1}$ for each $i \in \{1, 2, 3\}$. It is desired that the bound $a(\mathcal{F}) |V(G)|$ is attained by infinitely many non-isomorphic connected graphs, as is the case for $\mathcal{F} \in \{\mathcal{F}_{0,1}, \mathcal{F}_{1,k}, \mathcal{F}_{2,k}, \mathcal{F}_{3,k}\}$ by Theorem~\ref{mainresult1gen} (note that $\mathcal{F}_{0,1} \subseteq \mathcal{F}_{1,2} \cap \mathcal{F}_{2,2}$). Suppose that this holds for $\mathcal{F} = \mathcal{F}_{0,k}$ with $k \geq 2$. By any of Theorems~\ref{result3}--\ref{mainresult1gen}, $a(\mathcal{F}_{0,k}) \leq \frac{1}{k+1}$. In Section~\ref{Starsection}, we prove the following result. 

\begin{theorem}\label{startheorem} For $n \geq 0$ and $k \geq 2$, a connected $n$-vertex graph $G$ attaining the upper bound $\frac{n}{k+1}$ on $\iota(G, K_{1,k})$ exists if and only if $n = i(k+1)$ for some $i \in \{0, 1, 2\}$.
\end{theorem}
By Theorem~\ref{startheorem}, we actually have $a(\mathcal{F}_{0,k}) < \frac{1}{k+1}$. In Section~\ref{Starsection}, we also prove that $a(\mathcal{F}_{0,k}) \geq \frac{1}{k+\frac{3}{2}}$ (see Lemma~\ref{grapheq3}). 

In Section~\ref{Problemsection}, we pose some problems and conjectures, and also prove some results, motivated by the above. In particular, we establish additional intriguing phenomena concerning $k$-star isolation and $k$-cycle isolation.

\section{Proofs of Theorems~\ref{mainresult1} and \ref{mainresult1gen}} \label{Proofsection}

We start the proof of Theorem~\ref{mainresult1} with a lemma from \cite{Borg}. 

\begin{lemma}[\cite{Borg}] \label{lemma}
If $G$ is a graph, $\mathcal{F}$ is a set of graphs, $X \subseteq V(G)$ and $Y \subseteq N[X]$, then \[\iota(G, \mathcal{F}) \leq |X| + \iota(G-Y, \mathcal{F}).\] 
\end{lemma}
\textbf{Proof.} Let $D$ be an $\mathcal{F}$-isolating set of $G-Y$ of size $\iota(G-Y, \mathcal{F})$. Clearly, $\emptyset \neq V(F) \cap Y \subseteq V(F) \cap N[X]$ for each $\mathcal{F}$-graph $F$ that is a subgraph of $G$ and not a subgraph of $G-Y$. Thus, $X \cup D$ is an $\mathcal{F}$-isolating set of $G$. The result follows.~\hfill{$\Box$}
\\

We will say that a set $\mathcal{F}$ of graphs is \emph{component-represented} if each member of $\mathcal{F}$ has at least one component that is an $\mathcal{F}$-graph. Note that the union of component-represented sets of graphs is component-represented. Clearly, $\mathcal{F}_{0,k}$, $\mathcal{F}_{1,k}$, $\mathcal{F}_{2,k}$ and $\mathcal{F}_{3,k}$ are component-represented. 

\begin{lemma} \label{lemmacomp} If $G_1, \dots, G_r$ are the distinct components of a graph $G$, and $\mathcal{F}$ is a compo-nent-represented set of graphs, then $\iota(G,\mathcal{F}) = \sum_{i=1}^r \iota(G_i,\mathcal{F})$.
\end{lemma}
\textbf{Proof.} For each $i \in [r]$, let $D_i$ be a smallest $\mathcal{F}$-isolating set of $G_i$. Consider any $\mathcal{F}$-graph $F$ contained by $G$. Then, $F$ has a component $H$ that is an $\mathcal{F}$-graph. Since $H$ is connected, there exists some $j \in [r]$ such that $G_j$ contains $H$, so $N[D_j] \cap V(H) \neq \emptyset$, and hence $N[D_j] \cap V(F) \neq \emptyset$. Thus, $\bigcup_{i = 1}^r D_i$ is an $\mathcal{F}$-isolating set of $G$, and hence $\iota(G, \mathcal{F}) \leq \sum_{i = 1}^r |D_i| = \sum_{i = 1}^r \iota(G_i, \mathcal{F})$. Let $D$ be a smallest $\mathcal{F}$-isolating set of $G$. For each $i \in [r]$, $D \cap V(G_i)$ is an $\mathcal{F}$-isolating set of $G_i$. We have $\sum_{i = 1}^r \iota(G_i, \mathcal{F}) \leq \sum_{i = 1}^r |D \cap V(G_i)| = |D| = \iota(G, \mathcal{F})$. The result follows.~\hfill{$\Box$}
\\

The equation in Lemma~\ref{lemmacomp} may not hold if $\mathcal{F}$ is not component-represented. For example, if $G$ is the union of two vertex-disjoint $4$-paths $G_1$ and $G_2$ (the components of $G$), $F = (\{1, 2\}, \emptyset)$ and $\mathcal{F} = \{F\}$, then $\iota(G, \mathcal{F}) = 3$ and $\iota(G_1, \mathcal{F}) = \iota(G_2, \mathcal{F}) = 1$. We mention in passing that the following immediate consequence of Lemma~\ref{lemmacomp} has been used extensively in the literature because almost all the sets $\mathcal{F}$ that have been treated so far consist of connected graphs.

\begin{lemma} \label{lemmacompcor} If $G_1, \dots, G_r$ are the distinct components of a graph $G$, and $\mathcal{F}$ is a set of connected graphs, then $\iota(G,\mathcal{F}) = \sum_{i=1}^r \iota(G_i,\mathcal{F})$.
\end{lemma}

\noindent
\textbf{Proof of Theorem~\ref{mainresult1}.} Let $\mathcal{F} = \mathcal{F}_{0,k} \cup \mathcal{F}_{1,k}$. We first settle the second part of the theorem. The case $n \leq k-1$ is trivial. If $G = B_{n,k}$, then $n \neq k$ as $B_{k,k}$ is a $k$-clique and $(G,k)$ is not special. Suppose $n \geq k+1$. Let $B_1, \dots, B_q$ be the $(k+1)$-cliques in Construction~\ref{const1}. For each $j \in [q]$, let $B_j' = B_j - b_j$. Then, $B_1', \dots, B_q'$ are $(k-1)$-regular. For each $v \in V(B_{n,k})$, $N[v]$ does not intersect more than one of the vertex sets of $B_1', \dots, B_q'$,  so $\iota(B_{n,k}, \mathcal{F}) \geq q$. Since $\{b_1, \dots, b_q\}$ is an $\mathcal{F}$-isolating set of $B_{n,k}$, $\iota(B_{n,k}, \mathcal{F}) = q$.

Using induction on $n$, we now prove that the bound in the theorem holds. Since $\iota(G, \mathcal{F})$ is an integer, it suffices to prove that $\iota(G, \mathcal{F}) \leq \frac{n}{k+1}$. If $k \leq 2$, then a $K_k$-isolating set of $G$ is an $\mathcal{F}$-isolating set of $G$, so the result is given by Theorem~\ref{result1}. Consider $k \geq 3$. The result is trivial if $n \leq 2$ or $\iota(G,\mathcal{F}) = 0$. Suppose $n \geq 3$ and $\iota(G,\mathcal{F}) \geq 1$. 

Suppose $\Delta(G) \leq k-1$. Thus, $G$ contains no $k$-star. Since $\iota(G,\mathcal{F}) \geq 1$, $G$ contains a $(k-1)$-regular graph $R$. For each $v \in V(R)$, we have $k-1 = d_R(v) \leq d_G(v) \leq k-1$, so $N_G(v) = N_R(v)$. Thus, $E(G) \cap {V(R) \choose 2} = E(R)$. Suppose $V(G) \backslash V(R) \neq \emptyset$. Since $G$ is connected, we obtain that $vw \in E(G)$ for some $v \in V(R)$ and some $w \in V(G) \backslash V(R)$, which contradicts $N_G(v) = N_R(v)$. Thus, $V(G) \backslash V(R) = \emptyset$, which immediately yields $G = R$. Let $v \in V(G)$. Let $G' = G - N[v]$. Since $G$ is $(k-1)$-regular and is not a $k$-clique (as $(G,k)$ is not special), we have $n \geq k+1$, so $V(G') \neq \emptyset$. Suppose that $G'$ contains a $(k-1)$-regular graph $R'$. Since $G$ is connected, $uw \in E(G)$ for some $u \in V(R')$ and some $w \in V(G) \backslash V(R')$. We have $d_G(u) \geq d_{R'}(u) + 1 \geq k$, which contradicts $\Delta(G) \leq k-1$. Thus, $G'$ contains no $(k-1)$-regular graph. Since $\Delta(G') \leq \Delta(G) \leq k-1$, $G'$ contains no $k$-star. Thus, $\iota(G, \mathcal{F}) = 1 \leq \frac{n}{k+1}$. 

Now suppose $\Delta(G) \geq k$. Let $v \in V(G)$ with $d(v) = \Delta(G)$. If $V(G) = N[v]$, then $\{v\}$ is an $\mathcal{F}$-isolating set of $G$, so $\iota(G,\mathcal{F}) = 1 \leq \frac{n}{k+1}$. Suppose $V(G) \neq N[v]$. Let $G' = G-N[v]$ and $n' = |V(G')|$. Then, 
\[n \geq n' + k + 1\] 
and $V(G') \neq \emptyset$. Let $\mathcal{H}$ be the set of components of $G'$. Let $\mathcal{H}' = \{H \in \mathcal{H} \colon H \simeq K_k\}$. By the induction hypothesis, $\iota(H,\mathcal{F}) \leq \frac{|V(H)|}{k+1}$ for each $H \in \mathcal{H} \backslash \mathcal{H}'$ (recall that $k \geq 3$, meaning that $(H, k)$ is not special). If $\mathcal{H}' = \emptyset$, then by Lemma~\ref{lemma} (with $X = \{v\}$ and $Y = N[v]$) and Lemma~\ref{lemmacomp},
\begin{equation}
\iota(G,\mathcal{F}) \leq 1 + \iota(G',\mathcal{F}) = 1 + \sum_{H \in \mathcal{H}} \iota(H,\mathcal{F}) \leq 1 + \sum_{H \in \mathcal{H}} \frac{|V(H)|}{k+1} = \frac{k+1 + n'}{k+1} \leq \frac{n}{k+1}. \nonumber
\end{equation}

Suppose $\mathcal{H}' \neq \emptyset$. For any $H \in \mathcal{H}$ and any $x \in N(v)$ such that $xy_{x,H} \in E(G)$ for some $y_{x,H} \in V(H)$, we say that $H$ is \emph{linked to $x$} and that $x$ is \emph{linked to $H$}. Since $G$ is connected, each member of $\mathcal{H}$ is linked to at least one member of $N(v)$. For each $x \in N(v)$, let $\mathcal{H}'_x = \{H \in \mathcal{H}' \colon H \mbox{ is linked to } x\}$ and $\mathcal{H}_x^* = \{H \in \mathcal{H} \backslash \mathcal{H}' \colon H \mbox{ is linked to $x$ only}\}$. For each $H \in \mathcal{H} \backslash \mathcal{H}'$, let $D_H$ be an $\mathcal{F}$-isolating set of $H$ of size $\iota(H,\mathcal{F})$.
\\

\noindent
\emph{Case 1: $|\mathcal{H}'_x| \geq 2$ for some $x \in N(v)$.} For each $H \in \mathcal{H}' \backslash \mathcal{H}'_x$, let $x_H \in N(v)$ such that $H$ is linked to $x_H$. Let $X = \{x_H \colon H \in \mathcal{H}' \backslash \mathcal{H}'_x\}$. Note that $x \notin X$. Let 
\[D = \{v, x\} \cup X \cup \bigcup_{H \in \mathcal{H} \backslash \mathcal{H}'} D_{H} .\]
We have $V(G) = N[v] \cup \bigcup_{H \in \mathcal{H}} V(H)$, $y_{x,H} \in N[x]$ for each $H \in \mathcal{H}'_x$, and $y_{x_H,H} \in N[x_H]$ for each $H \in \mathcal{H}' \backslash \mathcal{H}'_x$, so $D$ is an $\mathcal{F}$-isolating set of $G$. Since $\iota(G,\mathcal{F}) \leq |D|$ and
\begin{align} n &= |N[v]| + k|\mathcal{H}'_x| + k|\mathcal{H}' \backslash \mathcal{H}_x'| + \sum_{H \in \mathcal{H} \backslash \mathcal{H}'} |V(H)| \nonumber \\
&\geq |\{v,x\} \cup X| + 2k + k|X| + \sum_{H \in \mathcal{H} \backslash \mathcal{H}'} (k+1)\iota(H,\mathcal{F}) \nonumber \\
&= 2(k+1) + (k+1)|X| + \sum_{H \in \mathcal{H} \backslash \mathcal{H}'} (k+1)|D_H| = (k+1)|D|, \nonumber
\end{align}
$\iota(G,\mathcal{F}) \leq \frac{n}{k+1}$.\\

\noindent
\emph{Case 2:}
\begin{equation} |\mathcal{H}'_x| \leq 1 \mbox{ \emph{for each} } x \in N(v). \label{k=3Hx<2} 
\end{equation} 
Let $H \in \mathcal{H}'$. Let $x \in N(v)$ such that $H$ is linked to $x$. Let $y = y_{x,H}$ and $X = \{x\} \cup V(H)$. Let $G^* = G - X$. Then, $G^*$ has a component $G_v^*$ such that $N[v] \backslash \{x\} \subseteq V(G_v^*)$, and the other components of $G^*$ are the members of $\mathcal{H}_{x}^*$. Let $D^*$ be an $\mathcal{F}$-isolating set of $G_v^*$ of size $\iota(G_v^*,\mathcal{F})$. Let $D = D^* \cup \{y\} \cup \bigcup_{I \in \mathcal{H}_{x}^*} D_I$. By Lemma~\ref{lemma}, since $X \subseteq N[y]$, $D$ is an $\mathcal{F}$-isolating set of $G$. We have
\begin{equation} \iota(G,\mathcal{F}) \leq |D^*| + 1 + \sum_{I \in \mathcal{H}_{x}^*} |D_I| \leq \iota(G_v^*,\mathcal{F}) + \frac{|X|}{k+1} + \sum_{I \in \mathcal{H}_{x}^*} \frac{|V(I)|}{k+1}.  \nonumber
\end{equation}
This yields $\iota(G,\mathcal{F}) \leq \frac{n}{k+1}$ if $\iota(G_v^*,\mathcal{F}) \leq \frac{|V(G_v^*)|}{k+1}$. Suppose $\iota(G_v^*,\mathcal{F}) > \frac{|V(G_v^*)|}{k+1}$. By the induction hypothesis, $G_v^* \simeq K_k$. Since $|N[v]| \geq k+1$ and $N[v] \backslash \{x\} \subseteq V(G^*_v)$, $V(G^*_v) = N[v] \backslash \{x\}$. Let $Y = (X \cup V(G^*_v)) \backslash \{v,x,y\}$ and $G_Y = G - \{v,x,y\}$. Then, the components of $G_Y$ are the components of $G[Y]$ and the members of $\mathcal{H}_x^*$.

Suppose that $G[Y]$ contains no $\mathcal{F}$-graph. Since $v, y \in N[x]$, $\{x\} \cup \bigcup_{I \in \mathcal{H}_x^*} D_I$ is an $\mathcal{F}$-isolating set of $G$, so
\begin{equation}
\iota(G,\mathcal{F}) \leq 1 + \sum_{I \in \mathcal{H}_x^*} |D_I| < \frac{|N[v]| + |V(H)|}{k+1} + \sum_{I \in \mathcal{H}_x^*} \frac{|V(I)|}{k+1} = \frac{n}{k+1}. \nonumber
\end{equation}

Now suppose that $G[Y]$ contains an $\mathcal{F}$-graph. Then, $G[Y]$ has a subgraph $F_Y$ that is a $k$-star or a $(k-1)$-regular graph (note that $G[Y]$ contains a $k$-star if it contains a regular graph of degree at least $k$). Thus, $|N_{G[Y]}[z]| \geq k$ for some $z \in V(F_Y)$. Let $W \subseteq N_{G[Y]}[z]$ such that $z \in W$ and $|W| = k$. Let $G_1 = G_v^*$, $G_2 = H$, $v_1 = v$, $v_2 = y$, $G_1' = G_1 - v_1$ and $G_2' = G_2 - v_2$. We have
\begin{equation}\label{edgesmain.2}
N_{G[Y]}[z] \subseteq Y = V(G_1') \cup V(G_2').
\end{equation}
Thus, $z \in V(G_j')$ for some $j \in \{1, 2\}$, and since $|V(G_1')| = |V(G_2')| = k-1 = |W|-1$, we have $|W \cap V(G_1')| \geq 1$ and $|W \cap V(G_2')| \geq 1$. Let $Z = V(G_j) \cup W$. Since $z$ is a vertex of the $k$-clique $G_j$,
\begin{equation} Z \subseteq N[z]. \label{edgesmain.3}
\end{equation}
We have
\begin{equation} |Z| = |V(G_j)| + |W \backslash V(G_j)| = k + |W \cap V(G_{3-j}')| \geq k+1. \label{edgesmain.4}
\end{equation}
Let $G_Z = G - Z$. Then, $V(G_Z) = \{x\} \cup (V(G_{3-j}) \backslash W) \cup \bigcup_{I \in \mathcal{H}_x^*} V(I)$. The components of $G_Z - x$ are $G_Z[V(G_{3-j}) \backslash W]$ (a clique having less than $k$ vertices) and the members of $\mathcal{H}_x^*$. Moreover, $v_{3-j} \in V(G_{3-j}) \backslash W$ (by (\ref{edgesmain.2})), $v_{3-j} \in N_{G_Z}(x)$, and by the definition of $\mathcal{H}_x^*$, $N_{G_Z}(x) \cap V(I) \neq \emptyset$ for each $I \in \mathcal{H}_x^*$. Thus, $G_Z$ is connected.
\\

\noindent
\emph{Subcase 2.1: $\mathcal{H}_x^* \neq \emptyset$.} Then, $G_Z$ is not a $k$-clique. By the induction hypothesis, $\iota(G_Z,\mathcal{F}) \leq \frac{|V(G_Z)|}{k+1}$. By (\ref{edgesmain.3}) and Lemma~\ref{lemma}, $\iota(G,\mathcal{F}) \leq 1 + \iota(G_Z,\mathcal{F}) \leq 1 + \frac{|V(G_Z)|}{k+1}$. By (\ref{edgesmain.4}), $\iota(G,\mathcal{F}) \leq \frac{|Z|}{k+1} + \frac{|V(G_Z)|}{k+1} = \frac{n}{k+1}$.
\\

\noindent
\emph{Subcase 2.2: $\mathcal{H}_x^* = \emptyset$.} Then, $G^* = G^*_v$, so $V(G) = V(G^*_v) \cup \{x\} \cup V(H)$ and $n = 2k+1$. We have $\Delta(G) = d(v) = k$. Thus, by (\ref{edgesmain.3}) and (\ref{edgesmain.4}), $N[z] = Z = V(G_j) \cup \{w\}$ for some $w \in V(G_{3-j}')$, and $V(G-N[z]) = \{x\} \cup V(G_{3-j} - w)$. If $G-N[z]$ contains no $\mathcal{F}$-graph, then $\iota(G,\mathcal{F}) = 1 < \frac{n}{k+1}$. Suppose that $G-N[z]$ contains an $\mathcal{F}$-graph $F$. Since $|V(G-N[z])| = k$, $G-N[z] = F \simeq K_k$. Since $\Delta(G) = k$, we have $N(x) = \{v_{j}\} \cup V(G_{3-j} - w)$ and, since $z \in V(G_j)$ and $w \in N[z] \cap V(G_{3-j})$, $N[w] = \{z\} \cup V(G_{3-j})$. Thus, $V(G-N[w]) = \{x\} \cup V(G_j - z)$. Since $|V(G-N[w])| = k \geq 3$ and $N[x] \cap V(G_j') = \emptyset$, $\{w\}$ is an $\mathcal{F}$-isolating set of $G$, so $\iota(G,\mathcal{F}) = 1 < \frac{n}{k+1}$.~\hfill{$\Box$}
\\
\\
\textbf{Proof of Theorem~\ref{mainresult1gen}.} Let $D$ be a smallest $(\mathcal{F}_{0,k} \cup \mathcal{F}_{1,k})$-isolating set of $G$. By Theorem~\ref{mainresult1}, $|D| \leq \left \lfloor \frac{n}{k+1} \right \rfloor$. Let $G' = G - N[D]$. Then, no subgraph of $G'$ is a $k$-star or a $(k-1)$-regular graph. Thus, $\Delta(G') \leq k-1$, $G'$ contains no $k$-cliques, and if $k = 3$, then $G'$ contains no cycles. By Brooks' Theorem, $\chi(G') \leq k-1$, so $\chi(H) \leq k-1$ for each subgraph $H$ of $G'$. Therefore, $D$ is an $\mathcal{F}_{3,k}$-isolating set of $G$, and hence if $i \in \{0, 1, 2\}$, then $D$ is also an $\mathcal{F}_{i,k}$-isolating set of $G$.

If $B_1, \dots, B_q$ are the $(k+1)$-cliques in Construction~\ref{const1}, then $B_1 - b_1, \dots, B_q - b_q$ are $(k-1)$-regular and $k$-chromatic. Thus, as in the proof of Theorem~\ref{mainresult1}, if $i \in \{1, 2, 3\}$ and $G = B_{n,k}$, then $\iota(G, \mathcal{F}_{i,k}) = \left \lfloor \frac{n}{k+1} \right \rfloor$.~\hfill{$\Box$}

\section{Isolation of stars} \label{Starsection}

We now address the problem, concerning the $K_{1,k}$-isolation number for $k \geq 2$, that is described in the last part of Section~\ref{Introsection}. We abbreviate $\iota(G,K_{1,k})$ to $\iota_k(G)$. By any of Theorems~\ref{result3}--\ref{mainresult1gen}, $\iota_k(G) \leq \frac{n}{k+1}$ for any connected graph $G$. The bound is attained if $G$ is a $(k+1)$-vertex graph containing a $k$-star. In this section, we first give an explicit construction of a connected $2(k+1)$-vertex graph $C(k)$ that also attains the bound, hence verifying the sufficiency condition in Theorem~\ref{startheorem}, we then address the claim at the end of Section~\ref{Introsection} that $a(\{K_{1,k}\}) \geq \frac{1}{k+\frac{3}{2}}$ by providing an explicit construction, based on $C(k)$, of a connected $n$-vertex graph $B_{n,C(k)}$ such that $\iota_k(B_{n,C(k)}) = \left \lfloor \frac{2n}{2k+3} \right \rfloor$ for any $n \geq 2k+3$, and we finally prove Theorem~\ref{startheorem}.

Let mod$^*$ be the usual modulo operation with the exception that for any integers $m$ and $n \neq 0$, $mn \mbox{ mod$^*$ } n$ is $n$ instead of $0$. For $1 \leq r < n$, let $C_n^r$ be the graph with $V(C_n^r) = [n]$ and 
\[E(C_n^r) = \bigcup_{i=1}^{n} \left\{ \{i, (i+j) \mbox{ mod$^*$ } n\} \colon j \in \left[ r \right] \right\},\]
%
that is, the \emph{$r$\textsuperscript{th} power of $C_n$} (the graph with vertex set $V(C_n)$ and where, for every two distinct vertices $v$ and $w$, $v$ and $w$ are neighbours if and only if the distance between them in $C_n$ is at most $r$). 

\begin{construction}\label{const2} \emph{Consider any integer $k \geq 2$. If $k$ is even, then let $C(k) = C_{2k+2}^{k/2}$. If $k$ is odd, then let $C(k)$ be the graph with $V(C(k)) = [2k+2]$ and}
\begin{center} $E(C(k)) = E \left(C_{2k+2}^{(k-1)/2}\right) \cup \left\{ \left\{i,i + \frac{k+1}{2} \right\} \colon i \in \left[\frac{k+1}{2} \right] \cup \left( \left[ k+1 + \frac{k+1}{2} \right] \backslash [k+1] \right) \right\}.$
\end{center}
\end{construction}  

\begin{lemma} \label{grapheq2} For $k \geq 2$,
\begin{equation} \iota_k(C(k)) = \gamma(C(k)) = 2. \nonumber 
\end{equation}
\end{lemma}
\textbf{Proof.} Let $s = 2k+2$. If $k$ is odd, then let $I = \left[\frac{k+1}{2} \right] \cup \left( \left[ k+1 + \frac{k+1}{2} \right] \backslash [k+1] \right)$. For $i \in [s]$, $N_{C(k)}(i) = \{x_i, (x_i + 1) \mbox{ mod$^*$ } s, \dots, (x_i + k) \mbox{ mod$^*$ } s\}$, where $x_i = \left(i - \frac{k}{2} \right) \mbox{ mod$^*$ } s$ if $k$ is even, $x_i = \left(i - \frac{k-1}{2} \right) \mbox{ mod$^*$ } s$ if $k$ is odd and $i \in I$, and $x_i = \left(i - \frac{k+1}{2} \right) \mbox{ mod$^*$ } s$ if $k$ is odd and $i \notin I$. Taking $j_i = (i+k+1) \mbox{ mod$^*$ } s$, we therefore have $d_{C(k)-N_{C(k)}[i]}(j_i) = k$ and $N_{C(k)}[\{i, j_i\}] = V(C(k))$. Thus, $\{i\}$ is not a $K_{1,k}$-isolating set of $C(k)$, and $\{i, j_i\}$ is a dominating set of $C(k)$.~\hfill{$\Box$}

\begin{construction}\label{const3} \emph{Consider any $n, k \in \mathbb{N}$ with $k \geq 2$ and $n \geq 2k+3$. Let $q = \left \lfloor \frac{n}{2k+3} \right \rfloor$. Thus, $n = q(2k+3) + r$, where $0 \leq r \leq 2k+2$. Let $u_1, \dots, u_{q+r}$ be the vertices $1, \dots, q+r$ of $P_{q+r}$, respectively. Let $R = (\emptyset, \emptyset)$ if $0 \leq r \leq 1$, and let $R = (\{u_{q+j} \colon j \in [r]\}, \{u_{q+r}u_{q+j} \colon j \in [r-1]\})$ if $r \geq 2$. Thus, if $r \geq 2$, then $R \simeq K_{1,r-1}$ and $V(R) \cap V(P_{q+1}) = \{u_{q+1}\}$. Since $n \geq 2k+3$, $q \geq 1$.  Let $G_1, \dots, G_q$ be copies of $C(k)$ such that $P_{q+r}, G_1, \dots, G_q$ are pairwise vertex-disjoint. For each $i \in [q]$, let $v_{i,1}, \dots, v_{i,2k+2}$ be the vertices of $G_i$ corresponding to the vertices $1, \dots, 2k+2$ of $C(k)$, respectively. Let $B_{n,C(k)}$ be the connected $n$-vertex graph with $V(B_{n,C(k)}) = V(P_{q+r}) \cup \bigcup_{i=1}^q V(G_i)$ and $E(B_{n,C(k)}) = \{u_iv_{i,1} \colon i \in [q]\} \cup E(P_{q+t}) \cup E(R) \cup \bigcup_{i=1}^q E(G_i)$, where $t = \min\{1, r\}$.}
\end{construction}  

\begin{lemma} \label{grapheq3} For $k \geq 2$ and $n \geq 2k+3$,
\[\iota_k(B_{n,C(k)}) = \left \lfloor \frac{2n}{2k+3} \right \rfloor.\]
\end{lemma}
\textbf{Proof.} Consider Construction~\ref{const3}. Let $G = B_{n,C(k)}$. Let $D$ be a smallest $K_{1,k}$-isolating set of $G$. For $i \in [q]$, let $D_i = D \cap (\{u_i\} \cup V(G_i))$. Let $D_i' = D_i$ if $u_i \notin D_i$, and let $D_i' = (D_i \backslash \{u_i\}) \cup \{v_{i,1}\}$ if $u_i \in D_i$. Thus, $N[D_i] \cap V(G_i) \subseteq N[D_i'] \cap V(G_i)$. For each $v \in V(G) \backslash (\{u_i\} \cup V(G_i))$, $N[v] \cap V(G_i) = \emptyset$. Thus, $D_i'$ is a $K_{1,k}$-isolating set of $G_i$. By Lemma~\ref{grapheq2}, $2 \leq |D_i'| \leq |D_i|$. Let $D_R = D \cap \{u_{q+j} \colon j \in [r]\}$.  We have $\iota_k(G) = |D| = |D_R| + \sum_{i=1}^q |D_i| \geq |D_R| + 2q$. Let $X = \{u_1, \dots, u_q, v_{1,k+2}, \dots, v_{q,k+2}\}$. If $r \leq k+1$, then $X$ is a $K_{1,k}$-isolating set of $G$, so $\iota_k(G) = 2q = \frac{2(n-r)}{2k+3} = \left \lfloor \frac{2(n-r)}{2k+3} + \frac{2r}{2k+3} \right \rfloor = \left \lfloor \frac{2n}{2k+3} \right \rfloor$. Suppose $r \geq k+2$. Then, $1 < \frac{2r}{2k+3} < 2$ and $R - u_{q+1}$ contains a $k$-star. Since $N[v] \cap V(R - u_{q+1}) = \emptyset$ for each $v \in V(G) \backslash V(R)$, we obtain $D_R \neq \emptyset$, so $\iota_k(G) \geq 2q + 1$. Since $X \cup \{u_{q+r}\}$ is a $K_{1,k}$-isolating set of $G$, $\iota_k(G) = 2q + 1 = \left \lfloor \frac{2(n-r)}{2k+3} + \frac{2r}{2k+3} \right \rfloor = \left \lfloor \frac{2n}{2k+3} \right \rfloor$.~\hfill{$\Box$}
\\

Lemma~\ref{grapheq3} yields Proposition~\ref{propeq3}. We now prove Theorem~\ref{startheorem}, using Lemma~\ref{grapheq2} and the next lemma.

\begin{lemma} \label{isodom} If $G$ is a connected $n$-vertex graph with $\iota_k(G) = \frac{n}{k+1} \in \{0, 1, 2\}$, then $\iota_k(G) = \gamma(G)$.
\end{lemma}
\textbf{Proof.} Let $i = \frac{n}{k+1}$. The result is trivial if $i = 0$. Suppose $i \geq 1$. Then, $G$ contains a $k$-star, so $\Delta(G) \geq k$. Let $v \in V(G)$ with $d(v) = \Delta(G)$. If $i = 1$, then $n = k+1$, $N[v] = V(G)$, and hence $\gamma(G) = 1$. Suppose $i = 2$. Then, $G - N[v]$ contains a $k$-star, so $|N_{G - N[v]}[w]| \geq k+1$ for some $w\in V(G - N[v])$. Since $2(k+1) = n \geq |N[v]| + |N_{G - N[v]}[w]| \geq 2(k+1)$, we have $|N[v]| = |N_{G - N[v]}[w]| = k+1 = |V(G - N[v])|$ and $\Delta(G) = k$. Thus, $N[w] = V(G - N[v])$, and hence $\{v, w\}$ is a dominating set of $G$. For any $u \in V(G)$, $|V(G - N[u])| \geq n - \Delta(G) - 1 = k+1$, so $\gamma(G) > 1$. Therefore, $\gamma(G) = 2$. \hfill{$\Box$} 
\\

If $S$ is a $k$-star, $k \geq 2$ and $v$ is the vertex of $S$ such that $N_S[v] = V(S)$, then $v$ is called the \emph{center of $S$}. 
If a graph $G$ contains a $k$-star $S$, then $S$ is called a \emph{$k$-star of $G$}. If $X, Y \subseteq V(G)$, then the set $\{xy \in E(G) \colon x \in X, \, y \in Y\}$ is denoted by $E_G(X,Y)$. Where no confusion arises, we may abbreviate $E_G(X,Y)$ to $E(X,Y)$.
\\
\\
\textbf{Proof of Theorem~\ref{startheorem}.} We use induction on $n$. If $n = 0$, then $\iota_k(G) = 0 = \frac{n}{k+1}$. If $1 \leq n \leq k$, then $G$ contains no $k$-star, so $\iota_k(G) = 0 < \frac{n}{k+1}$. If $n = k+1$, then $\iota_k(G) \leq 1 = \frac{n}{k+1}$, and equality holds if $G$ contains a $k$-star. Suppose $n \geq k+2$. If $\Delta(G) \leq k-1$, then $G$ contains no $k$-star, so $\iota_k(G) = 0 < \frac{n}{k+1}$. Suppose $\Delta(G) \geq k$. Let $v_0 \in V(G)$ such that $d(v_0) = \Delta(G)$. Let $G' = G-N[v_0]$, and let $n' = |V(G')|$. If $G'$ contains no $k$-star, then $\iota_k(G) = 1 < \frac{n}{k+1}$. Suppose that $G'$ contains a $k$-star $S_1$. Then, $n \geq |N[v_0]| + |V(S_1)| \geq 2(k+1)$. Let $v_1$ be the center of $S_1$. If $n = 2(k+1)$, then $\{v_0, v_1\}$ is a dominating set of $G$, so $\iota_k(G) \leq 2 = \frac{n}{k+1}$, and by Lemma~\ref{grapheq2}, equality holds if $G = C(k)$. Suppose $n > 2(k+1)$. Let $G_1, \dots, G_r$ be the distinct components of $G'$, where $|V(G_1)| \geq \cdots \geq |V(G_r)|$. We have $n' = n - d(v_0) - 1 = \sum_{j=1}^r |V(G_j)|$. Since $G$ is connected, for each $j \in [r]$, 
\begin{equation} 
E(N(v_0), V(G_j)) \neq \emptyset. \label{G_jlinks}
\end{equation}

By the induction hypothesis, for each $j \in [r]$, $\iota_k(G_j) \leq\frac{|V(G_j)|}{k+1}$, and equality holds only if $|V(G_j)| \in \{k+1, 2(k+1)\}$. By Lemma~\ref{lemma} (with $X = \{v_0\}$ and $Y = N[v_0]$) and Lemma~\ref{lemmacomp}, $\iota_k(G) \leq 1 + \sum_{j=1}^r \iota_k(G_j)$. Thus, if $\Delta(G) > k$ or $\iota_k(G_{j'}) < \frac{|V(G_{j'})|}{k+1}$ for some $j' \in [r]$, then 
\[\iota_k(G) < \frac{|N[v_0]|}{k+1} + \sum_{j=1}^r \frac{|V(G_j)|}{k+1} = \frac{n}{k+1}.\]
Suppose $\Delta(G) = k$ and $\iota_k(G_j) = \frac{|V(G_j)|}{k+1}$ for each $j \in [r]$. Then, for each $j \in [r]$, $|V(G_j)| \in \{k+1, 2(k+1)\}$, and by Lemma~\ref{isodom}, $G_j$ has a dominating set $D_j$ of size $\frac{|V(G_j)|}{k+1}$.

Suppose $n > 3(k+1)$. Then, $r \geq 2$. Let $G^* = G - V(G_r)$. By (\ref{G_jlinks}), $G^*$ is connected. If $|V(G_r)| = k+1$, then $|V(G^*)| > 2(k+1)$. If $|V(G_r)| = 2(k+1)$, then, since $|V(G_1)| \geq |V(G_r)|$, $|V(G^*)| \geq 3(k+1)$. By the induction hypothesis, $G^*$ has a $K_{1,k}$-isolating set $D^*$ such that $|D^*| < \frac{|V(G^*)|}{k+1}$. Since $D^* \cup D_r$ is a $K_{1,k}$-isolating set of $G$, $\iota_k(G) < \frac{|V(G^*)|}{k+1} + \frac{|V(G_r)|}{k+1} = \frac{n}{k+1}$.

Now suppose $n \leq 3(k+1)$. Since $d(v_0) = \Delta(G) = k$ and $n' = \sum_{j=1}^r |V(G_j)| = h(k+1)$ for some integer $h \geq 1$, we have $n - (k + 1) = h(k+1)$, so $n = (h+1)(k+1)$. Since $2(k+1) < n \leq 3(k+1)$, we obtain $n = 3(k+1)$. Let $S_0$ be the $k$-star of $G$ with $V(S_0) = N[v_0]$ and $E(S_0) = \{v_0x \colon x \in N(v_0)\}$. Then, $G' = G - V(S_0)$. Let $G'' = G' - V(S_1)$. We have $|V(G'')| = n - |V(S_0)| - |V(S_1)| = 3(k+1) - 2(k+1) = k+1$. If $G''$ contains no $k$-star, then, since $G'' = G - N[\{v_0, v_1\}]$, $\{v_0, v_1\}$ is a $K_{1,k}$-isolating set of $G$, so $\iota_k(G) \leq 2 < \frac{n}{k+1}$. Suppose that $G''$ contains a $k$-star $S_2$. Then, $V(G'') = V(S_2)$, and $V(S_0)$, $V(S_1)$ and $V(S_2)$ form a partition of $V(G)$, that is,
\begin{equation} V(G) = V(S_0) \cup V(S_1) \cup V(S_2), \label{Gpartition}
\end{equation}
and $V(S_0)$, $V(S_1)$ and $V(S_2)$ are pairwise disjoint.

Let $X_p = V(S_p)$ for each $p \in \{0, 1, 2\}$. Since $G$ is connected, $E(X_p, X_q) \neq \emptyset$ for some $p, q \in \{0,1,2\}$ with $p \neq q$. Let $s$ be the unique member of $\{0, 1, 2\} \backslash \{p, q\}$. Since $G$ is connected, $E(X_s, X_t) \neq \emptyset$ for some $t \in \{p,q\}$. We may assume that $t = p = 0$. Thus, $E(X_0, X_1) \neq \emptyset \neq E(X_0, X_2)$, and hence $x_0x_1, x_0'x_2 \in E(G)$ for some $x_0, x_0' \in X_0$, $x_1 \in X_1$ and $x_2 \in X_2$. Since $N[v_0] = X_0$, we have $x_0 \neq v_0$ and $x_0' \neq v_0$, so $x_0, x_0' \in N(v_0)$. Similarly, $x_1 \in N(v_1)$ and $x_2 \in N(v_2)$, where $v_2$ is the center of $S_2$.

Recall that $\Delta(G) = k$. Let $M = \{v \in V(G) \colon d(v) = k\}$. We will prove the result by showing that there exists a subset $D$ of $V(G)$ such that 
\begin{equation} \mbox{$|D| \leq 2$ and $N[D] \cap N[v] \neq \emptyset$ for each $v \in M$.} \label{RTP} 
\end{equation}
This gives $\Delta(G-N[D]) \leq k-1$, so $G-N[D]$ contains no $k$-star, and hence $\iota_k(G) \leq 2 < \frac{n}{k+1}$, as required. We will often use the fact that if $u \in D \cap X_i$ for some $i \in \{0, 1, 2\}$, then 
\begin{equation} \mbox{$v_i \in N[D] \cap N[x]$ for each $x \in X_i$.} \label{commonfact}
\end{equation} 

For each $i \in \{0, 1, 2\}$, let $M^{(i)} = N(v_i) \cap M$. By (\ref{Gpartition}), 
\begin{equation} M = \{v_0, v_1, v_2\} \cup M^{(0)} \cup M^{(1)} \cup M^{(2)}. \label{Mpartition} 
\end{equation}
Let
\begin{align} M_0 &= \{x \in M^{(0)} \colon N[x] = X_0\}, \nonumber \\
M_1 &= \{x \in M^{(0)} \colon N(x) \cap X_1 \neq \emptyset, \, N(x) \cap X_2 = \emptyset\}, \nonumber \\
M_2 &= \{x \in M^{(0)} \colon N(x) \cap X_1 = \emptyset, \, N(x) \cap X_2 \neq \emptyset\}, \nonumber \\
M_3 &= \{x \in M^{(0)} \colon N(x) \cap X_1 \neq \emptyset, \, N(x) \cap X_2 \neq \emptyset\}. \nonumber
\end{align}
For each $x \in M^{(0)}$ with $N[x] \subseteq X_0$, we have $|N[x]| = k + 1 = |X_0|$, so $x \in M_0$. Thus, by (\ref{Gpartition}),
\begin{equation} M^{(0)} = M_0 \cup M_1 \cup M_2 \cup M_3.\bigskip \label{M(0)}
\end{equation}

\noindent\textit{Case 1: $M_1 = \emptyset$ or $M_2 = \emptyset$.} We may assume that $M_2 = \emptyset$, so $M^{(0)} = M_0 \cup M_1 \cup M_3$ by (\ref{M(0)}). Let $D = \{v_1, x_2\}$. Then, $x_0' \in N[D] \cap N[v]$ for each $v \in \{v_0\} \cup M_0$, and $\emptyset \neq N(v) \cap X_1 \subseteq N(v) \cap N[D]$ for each $v \in M_1 \cup M_3$. Together with (\ref{commonfact}) and (\ref{Mpartition}), this gives us that (\ref{RTP}) holds, as required.
\\

\noindent\textit{Case 2: $M_1 \neq \emptyset$ and $M_2 \neq \emptyset$.}
\\

\noindent\textit{Subcase 2.1: For some $i \in \{1, 2\}$, there exist some $y \in N(v_i)$ and $x \in M_{3-i}$ such that $N(y) \cap N(v_0) \nsubseteq N[x]$.} Recall that $x_0x_1 \in E(G)$. We may assume that 
\begin{equation} \mbox{$i = 1$, $y = x_1$ and $x_0 \in (N(y) \cap N(v_0)) \backslash N[x]$.} \label{case 2.1}
\end{equation}
Since $i = 1$, we have $x \in M_2$, so $N(x) \cap X_1 = \emptyset$ and $N(x) \cap X_2 \neq \emptyset$. We may assume that $x_2 \in N(x) \cap X_2$ (and $x = x_0'$).

Let $D_1 = \{x, x_1\}$. If (\ref{RTP}) holds with $D = D_1$, then we are done. Suppose that it does not. Then, by (\ref{commonfact}), $N[D_1] \cap N[u] = \emptyset$ for some $u \in M \cap X_2$. Since $x_2 \in N(x) \cap N(v_2)$, $u \in M^{(2)}$. 

Let $D_2 = \{x, u\}$. If (\ref{RTP}) holds with $D = D_2$, then we are done. Suppose that it does not. Then, by (\ref{commonfact}), $N[D_2] \cap N[w] = \emptyset$ for some $w \in M \cap X_1$. Let $I = V(G) \backslash N[D_2]$. Then, $N[w] \subseteq I$. Since $k+1 = |N[w]| \leq |I| = n - |N[x]| - |N[u]| = k+1$, $N[w] = I$. By (\ref{case 2.1}), $x_0 \notin N[x]$. Since $N[D_1] \cap N[u] = \emptyset$ and $x_0 \in N(x_1)$, $x_0 \notin N[u]$. Thus, $x_0 \notin N[D_2]$, and hence $x_0 \in N[w]$ (as $N[w] = I$).   

Let $D = \{x_0, u\}$. Since $x_0 \in N[w]$, $w \in N[D]$. Let $U = X_1 \cap N(u)$ and $U' = X_1 \backslash N(u)$. Since $N(x) \cap X_1 = \emptyset$ and $N[w] = I = V(G) \backslash (N[x] \cup N[u])$, we have $U' \subseteq N[w]$, so $w \in N[D] \cap N[v]$ for each $v \in U'$. For each $v \in U$, we have $v \in N(u)$, so $u \in N[D] \cap N[v]$. By (\ref{commonfact}), $N[D] \cap N[v] \neq \emptyset$ for each $v \in X_0 \cup X_2$. Thus, $N[D] \cap N[v] \neq \emptyset$ for each $v \in V(G)$, and hence (\ref{RTP}) holds.
\\

\noindent\textit{Subcase 2.2: For each $i \in \{1, 2\}$,} 
\begin{equation} \mbox{\textit{$N(y) \cap N(v_0) \subseteq N[x]$ for every $y \in N(v_i)$ and $x \in M_{3-i}$}.} \label{case 2.2}
\end{equation}
Recall that $M_1 \neq \emptyset$ and $M_2 \neq \emptyset$ (as we are in Case~2). Let $x \in M_2$. Thus, $N(x) \cap X_1 = \emptyset$ and $N(x) \cap X_2 \neq \emptyset$. We may assume that $x = x_0'$, giving $x_2 \in N(x) \cap X_2$. 

Let $D_1 = \{x, v_1\}$. If (\ref{RTP}) holds with $D = D_1$, then we are done. Suppose that it does not. Then, by (\ref{commonfact}), $N[D_1] \cap N[u] = \emptyset$ for some $u \in M \cap X_2$. Since $x_2 \in N(x) \cap N(v_2)$, $u \in M^{(2)}$. Let $I = V(G) \backslash N[D_1]$, $I_0 = X_0 \backslash N[x]$ and $I_2 = X_2 \backslash N[x]$. We have $N[u] \subseteq I$. Since $N[D_1] = N[x] \cup X_1$ and $k+1 = |N[u]| \leq |I| = n - |N[x]| - |N[v_1]| = k+1$, 
\begin{equation} N[u] = I = I_0 \cup I_2. \label{2.2fact1}
\end{equation}

Let $D_2 = \{x_1, v_2\}$. If (\ref{RTP}) holds with $D = D_2$, then we are done. Suppose that it does not. Then, by (\ref{commonfact}), $N[D_2] \cap N[w] = \emptyset$ for some $w \in M \cap X_0$. Since $x_0 \in N(x_1) \cap N(v_0)$, $w \in M^{(0)}$. 
We have $N[w] \subseteq V(G) \backslash N[D_2] = V(G) \backslash (N[x_1] \cup X_2) \subseteq (X_0 \backslash \{x_0\}) \cup X_1$, $|N[w]| = |X_0 \backslash \{x_0\}| + 1$, and hence $w \in M_1$. Thus, $x_1' \in N(w)$ for some $x_1' \in X_1$, and $I_0 \subseteq N[w]$ by (\ref{case 2.2}) and (\ref{2.2fact1}).

Let $D = \{x_1', x_2\}$. By (\ref{commonfact}), $N[D] \cap N[v] \neq \emptyset$ for each $v \in X_1 \cup X_2$. Since $I_0 \subseteq N[w]$, $w \in N[x_1'] \cap N[v]$ for each $v \in I_0$. Since $X_0 \backslash I_0 = N[x] \cap X_0$ (and $x_2 \in N(x)$), $x \in N[x_2] \cap N[v]$ for each $v \in X_0 \backslash I_0$. Thus, $N[D] \cap N[v] \neq \emptyset$ for each $v \in V(G)$, and hence (\ref{RTP}) holds.~\hfill{$\Box$}

\section{Problems and further results} \label{Problemsection}

We may assume that the vertex set of an $n$-vertex graph is $[n]$. Let 
\[\mathcal{G} = \{G \colon G \mbox{ is a connected graph, } V(G) = [n] \mbox{ for some } n \geq 1\}.\] 
Thus, $\mathcal{G}$ is an infinite set. For any set $\mathcal{F}$ of graphs and any real number $\alpha > 0$, let 
\[\mathcal{G}(\mathcal{F},\alpha) = \{G \in \mathcal{G} \colon \iota(G,\mathcal{F}) \leq \lfloor \alpha |V(G)| \rfloor\},\] 
and let
\[\mathcal{G}(\mathcal{F},\alpha)^* = \{G \in \mathcal{G}(\mathcal{F},\alpha) \colon \iota(G,\mathcal{F}) = \lfloor \alpha |V(G)| \rfloor\} \quad \mbox{and} \quad \mathcal{G}[\mathcal{F},\alpha] = \mathcal{G} \backslash \mathcal{G}(\mathcal{F},\alpha).\]
Thus, $\mathcal{G}[\mathcal{F},\alpha] = \{G \in \mathcal{G} \colon \iota(G,\mathcal{F}) > \lfloor \alpha |V(G)| \rfloor\}$. In view of Theorems~\ref{result1}--\ref{mainresult1gen} and the last part of Section~\ref{Introsection}, we pose the following problems.

\begin{problem} \label{problem1} (a) Is there a rational number $c(\mathcal{F})$ such that $\mathcal{G}{[\mathcal{F},c(\mathcal{F})]}$ is finite and $\mathcal{G}{(\mathcal{F},c(\mathcal{F}))}^*$ is infinite for \\
(i) $\mathcal{F} = \{K_{1,k}\}$? \\
(ii) $\mathcal{F} = \{C_k\}$? \\
(iii) $\mathcal{F} = \{P_k\}$?\medskip  
\\
(b) If $c(\mathcal{F})$ exists, then determine $c(\mathcal{F})$, $\mathcal{G}{[\mathcal{F},c(\mathcal{F})]}$ and (at least) an infinite subset of $\mathcal{G}(\mathcal{F},c(\mathcal{F}))^*$.
\end{problem}
By Theorem~\ref{mainresult1gen}, for each $i \in \{1, 2, 3\}$, $c(\mathcal{F}_{i,k}) = \frac{1}{k+1}$, $\mathcal{G}{[\mathcal{F}_{i,k},c(\mathcal{F}_{i,k})]} = \{G \in \mathcal{G} \colon (G,k)$ is special$\}$ and $\{G \in \mathcal{G} \colon G \simeq B_{n,k}, \, n \geq 3, n \neq k\} \subseteq \mathcal{G}{(\mathcal{F}_{i,k},c(\mathcal{F}_{i,k}))^*}$. 

\begin{conj} For each of (i)--(iii) of Problem~\ref{problem1}(a), $c(\mathcal{F})$ exists.
\end{conj}

We may abbreviate $c(\{F\})$ to $c(F)$. By Ore's result and the Caro--Hansberg--\.{Z}yli\'{n}ski result (the cases $k = 1$ and $k = 2$ of Theorem~\ref{result1}, respectively), 
\begin{equation} c(K_{1,k}) = \frac{1}{k+2} \quad \mbox{ for } 0 \leq k \leq 1. \label{constK1k}   
\end{equation}
Suppose that $c(K_{1,k})$ exists for $k \geq 2$. By any of Theorems~\ref{result3}--\ref{mainresult1gen}, $c(K_{1,k}) \leq \frac{1}{k+1}$. We now show that, rather surprisingly, $c(K_{1,k}) \geq \frac{1}{k + \frac{3}{2}} = \frac{2}{2k+3}$ for $k \geq 2$.

\begin{prop}\label{propeq3} If $k \geq 2$ and $c(K_{1,k})$ exists, then
\[c(K_{1,k}) \geq \frac{2}{2k+3}.\]
\end{prop}
\textbf{Proof.} Let $\alpha$ be a real number such that $0 < \alpha < \frac{2}{2k+3}$. For each $r \in \mathbb{N}$, $\lfloor \alpha r(2k+3)\rfloor \leq \alpha r(2k+3) < 2r = \iota(B_{r(2k+3),C(k)}, K_{1,k})$ by Lemma~\ref{grapheq3}. Thus, $\{B_{r(2k+3),C(k)} \colon r \in \mathbb{N}\}$ is an infinite subset of $\mathcal{G}[\{K_{1,k}\}, \alpha]$.~\hfill{$\Box$}

\begin{conj} For $k \geq 2$,
\[c(K_{1,k}) = \frac{2}{2k+3}.\]
\end{conj}

Recall that $C_k$ was defined for $k \geq 3$. Let $C_1$ and $C_2$ be $K_1$ and $K_2$, respectively. Then, $C_k = ( [k], \{ ij \in {[k] \choose 2} \colon j = (i+1) \mbox{ mod$^*$ } k \} )$ for $k \geq 1$. By Ore's result, the Caro--Hansberg--\.{Z}yli\'{n}ski result and Theorem~\ref{result2}, 
\begin{equation} c(C_k) = \frac{1}{k+1} \quad \mbox{for } 1 \leq k \leq 3. \label{constCk}  
\end{equation}
This is also given by Theorem~\ref{result1}. Bartolo, Scicluna and the present author \cite{BBS} recently showed that $c(C_k) = \frac{1}{k+1}$ also holds for $k = 4$. We now show that, surprisingly, if $k \geq 5$ and $c(C_k)$ exists, then $c(C_k) \geq \frac{1}{k + \frac{1}{2}} = \frac{2}{2k+1}$.

\begin{lemma} \label{grapheq4} For $k \geq 4$,
\begin{equation} \iota(C(k), C_{k+1}) = 2. \nonumber 
\end{equation}
\end{lemma}
\textbf{Proof.} Let $s$, $I$ and $j_1, \dots, j_s$ be as in the proof of Lemma~\ref{grapheq2}.  Consider any $i \in [s]$. Let $\ell = \lceil (k-1)/2 \rceil$ if either $k$ is even or $k$ is odd and $j_i \in I$, and let $\ell = (k+1)/2$ if $k$ is odd and $j_i \notin I$. Let $p = k - \ell$. Since $k \geq 4$, we have $\ell \geq 2$ and $p \geq 2$. For $r \in [\ell]$, let $x_r = (j_i - r) \mbox{ mod$^*$ } s$. For $r \in [p]$, let $y_r = (j_i + r) \mbox{ mod$^*$ } s$. Let $H = C(k) - N_{C(k)}[i]$. Note that $V(H) = N_{C(k)}[j_i]$ and that $H$ contains the $(k+1)$-cycle $(V(H), \{j_i x_{\ell}, x_{\ell} x_{\ell - 1}, \dots, x_2x_1, x_1y_1, y_1y_2, \dots, y_{p-1}y_p, y_p j_i\})$. Thus, $\{i\}$ is not a $C_{k+1}$-isolating set of $C(k)$, and $\{i, j_i\}$ is a dominating set of $C(k)$. \hfill{$\Box$} 

\begin{construction}\label{const4} \emph{For $k \geq 2$ and $n \geq 2k+3$, let $B_{n,C(k)}'$ be the graph with $V(B_{n,C(k)}') = V(B_{n,C(k)})$ and $E(B_{n,C(k)}') = E(B_{n,C(k)}) \cup {V(R) \choose 2}$, where $R$ is as in Construction~\ref{const3}.}
\end{construction}

\begin{lemma} \label{grapheq5} For $k \geq 4$ and $n \geq 2k+3$,
\[\iota(B_{n,C(k)}', C_{k+1}) = \left \lfloor \frac{2n}{2k+3} \right \rfloor.\]
\end{lemma}
\textbf{Proof.} The result is obtained by an argument similar to that in the proof of Lemma~\ref{grapheq3}, using Lemma~\ref{grapheq4}.~\hfill{$\Box$}

\begin{prop} If $k \geq 5$ and $c(C_k)$ exists, then
\[c(C_k) \geq \frac{2}{2k+1}.\]
\end{prop}
\textbf{Proof.} Let $\alpha$ be a real number such that $0 < \alpha < \frac{2}{2k+1}$. For each $r \in \mathbb{N}$, $\lfloor \alpha r(2k+1)\rfloor \leq \alpha r(2k+1) < 2r = \iota(B_{r(2k+1),C(k-1)}', C_k)$ by Lemma~\ref{grapheq5}. Thus, $\{B_{r(2k+1),C(k-1)}' \colon r \in \mathbb{N}\}$ is an infinite subset of $\mathcal{G}[\{C_k\}, \alpha]$.~\hfill{$\Box$}\\

Our next problem is stronger than Problem~\ref{problem1}(a).

\begin{problem} What is the smallest rational number $c(\mathcal{F},n)$ such that $\iota(G, \mathcal{F}) \leq c(\mathcal{F},n)n$ for every connected $n$-vertex graph $G$ if \\
(i)~$\mathcal{F} = \{K_{1,k}\}$? \\
(ii)~$\mathcal{F} = \{C_k\}$? \\ 
(iii)~$\mathcal{F} = \{P_k\}$?
\end{problem}
By Theorem~\ref{mainresult1gen}, for each $i \in \{1, 2, 3\}$, $c(\mathcal{F}_{i,k},k) = \frac{1}{k}$, $c(\mathcal{F}_{i,2},5) = \frac{2}{5}$, and $c(\mathcal{F}_{i,k},n) = \left \lfloor \frac{n}{k+1} \right \rfloor \frac{1}{n}$ for $(n,k) \notin \{(k,k), (5,2)\}$.

\begin{problem} \label{problem2} Determine $\mathcal{G}(\mathcal{F}_{i,k},c(\mathcal{F}_{i,k}))^*$ for $i \in \{1, 2, 3\}$.
\end{problem}

\noindent
\textbf{Acknowledgements.} The author is grateful to the anonymous referees for checking the paper and providing constructive remarks that led to an improvement in the presentation.

\end{document}